\documentclass[final]{elsarticle}  

\usepackage[colorlinks,citecolor=blue,linktoc=all,linkcolor=cyan]{hyperref}
\usepackage{graphicx}

\usepackage[T1]{fontenc}
\usepackage{dsfont}               
\usepackage{mathrsfs}             
\usepackage{slashed}              
\usepackage{amsmath}
\usepackage{amssymb}
\usepackage{amsbsy}
\usepackage{amsfonts}

\numberwithin{equation}{section}

\journal{arXiv}

\usepackage{titlesec}
\usepackage{sectsty}
\titleformat{\section}{\normalfont\Large\bfseries}{\thesection}{1em}{}
\titleformat{\subsection}{\normalfont\large\bfseries}{\thesubsection}{1em}{}
\titleformat{\subsubsection}{\normalfont\normalsize\bfseries}{\thesubsubsection}{1em}{}

\newtheorem{lema}{Lemma}
\newproof{prova}{Proof}
\newtheorem{teo}{Theorem}

\begin{document}

\begin{frontmatter}

\title{Beukers-like proofs of irrationality for $\zeta{(2)}$ and $\zeta{(3)}$}

\author[endereco]{F.~M.~S.~Lima}

\ead{fmsl@unb.br}

\address[endereco]{Institute of Physics, University of Bras\'{i}lia, 70919-970, Bras\'{i}lia-DF, Brazil}

\begin{abstract}
In this note, I develop step-by-step proofs of irrationality for $\,\zeta{(2)}\,$ and $\,\zeta{(3)}$. Though the proofs follow closely those based upon unit-square integrals proposed originally by Beukers, I introduce some modifications which certainly will be useful for those interested in understanding this kind of proof and/or trying to extend it to higher zeta values, Catalan's constant, or other related numbers.
\end{abstract}


\end{frontmatter}





\section{Introduction}

For real values of $\,s$, $s>1$, the Riemann zeta function is defined as $\,\zeta(s) := \sum_{n=1}^\infty{1/n^s}$.\footnote{In this domain, this series converges to a real number greater than $\,1/(s-1)$, according to the integral test. Indeed, since $\,\sum_{n=1}^\infty{1/n^s} = 1 +\sum_{n=2}^\infty{1/n^s}$, then $\,\zeta{(s)} > 1$.}  In 1978, R. Ap\'{e}ry succeeded in proofing that both $\,\zeta{(2)} = \sum_{n=1}^\infty{1/n^2}\,$ and $\,\zeta{(3)} = \sum_{n=1}^\infty{1/n^3}\,$ are irrational numbers~\cite{Apery,Poorten}.\footnote{It is well-known that $\,\zeta{(2)} = \pi^2/6$, as first proved by Euler (for distinct proofs, see Ref.~\cite{Lima2012} and references therein). Of course, this imply that $\,\zeta{(2)}\,$ is irrational, since $\,\pi\,$ is a transcendental number, as first proved by Lindemann (1882), but Ap\'{e}ry's proof is still interesting because it does not use any mathematical property of $\,\pi$, as well as because it gives an estimate for the irrationality measure of $\,\pi^2$ (we shall not explore this measure here).}  His proofs were soon shortened by F.~Beukers (1979)~\cite{Beukers79}, who translated it into equivalent statements about certain improper integrals over the unit square $[0,1]^2$. In this note, I develop rigorous Beukers-like proofs of irrationality for $\zeta{(2)}$ and $\zeta{(3)}$, in full details. Hopefully, this complete version of irrationality proofs can be useful for those interested in understanding this kind of proof and developing further generalizations.

\section{Irrationality of $\zeta{(2)}$}

\subsection{Preliminaries}

We begin with some non-trivial lemmas on certain \emph{improper} double integrals which are essential for the Beukers' irrationality proofs.

\begin{lema}[An unit square integral for $\,\zeta{(2)}\,$]
\label{lem:I0}
\begin{equation*}
\int_0^1\!\!\!\int_0^1{\frac{1}{1-x y} \: dx \, dy} = \zeta{(2)} \, .
\label{eq:I0}
\end{equation*}
\end{lema}

\begin{prova}
\; Let us define $\,I_{00} := \int_0^1\!\!\int_0^1{1/(1-x y) \: dx \, dy}$.  Since the integrand tends to infinity as $\,(x,y) \rightarrow (1,1)$, then $I_{00}$ is an \emph{improper} integral and we need to evaluate the corresponding limit, namely
\begin{equation}
I_{00} = \lim_{\varepsilon \rightarrow \, 0^{+}}\int_0^{1-\varepsilon}\!\!\int_0^{1-\varepsilon}{\!\frac{1}{1-x y} ~ dx \, dy} .
\end{equation}
By expanding the integrand as a geometric series, one finds
\begin{eqnarray}
I_{00} &=& \int_0^{1^{-}}\!\!\!\int_0^{1^{-}}{\!\sum_{n=0}^\infty{(x y)^n} \: dx \, dy}  \nonumber \\
&=& \sum_{n=0}^\infty{\int_0^{1^{-}}\!\!\!\int_0^{1^{-}}{\!x^n y^n \: dx \, dy}} = \sum_{n=0}^\infty{ \left( \int_0^{1^{-}}{\!x^n \: dx} \, \int_0^{1^{-}}{\!y^n \: dy}\right)}  \nonumber \\
&=& \sum_{n=0}^\infty{ \left\{\left[ \frac{x^{n+1}}{n+1} \right]_0^{1^{-}} \cdot \left[ \frac{y^{n+1}}{n+1} \right]_0^{1^{-}}\right\}} = \sum_{n=0}^\infty{ \, \lim_{\varepsilon \rightarrow \, 0^{+}}{ \left( \frac{(1-\varepsilon)^{n+1}}{n+1} \right)^{\!2}}}  \nonumber \\
&=& \sum_{n=0}^\infty{ \frac{1}{(n+1)^2} } = \sum_{m=1}^\infty{ \frac{1}{m^2}} = \zeta{(2)} \, .
\end{eqnarray}
The interchange of limits, series and integrals is fully justified, within the rigor of mathematical analysis, in the proof of Theorem 2.1 of Ref.~\cite{Hadji2002}.
\begin{flushright} $\Box$ \end{flushright}
\end{prova}

\begin{lema}[$I_{rr}$]
\label{lem:Ir}
\; For all integers $\,r>0$
\begin{equation*}
\int_0^1\!\!\!\int_0^1{x^r y^r \frac{1}{1-x y} \: dx \, dy} = \zeta{(2)} -\sum_{m=1}^r{\frac{1}{m^2}} \, .
\label{eq:Ir}
\end{equation*}
\end{lema}

\begin{prova}
\; As in the previous lemma, define $\,I_{rr} := \int_0^1\!\!\int_0^1{x^r y^r/(1-x y) \, dx\,dy}$.  It follows that
\begin{eqnarray}
I_{rr} &=& \lim_{\varepsilon \rightarrow \, 0^{+}}\int_0^{1-\varepsilon}\!\!\int_0^{1-\varepsilon}{\frac{x^r y^r}{1-x y} \: dx \, dy} \nonumber \\
&=& \int_0^{1^{-}}\!\!\!\int_0^{1^{-}}{x^r y^r \, \sum_{n=0}^\infty{(x y)^n} \: dx \, dy} = \sum_{n=0}^\infty{\int_0^{1^{-}}\!\!\!\int_0^{1^{-}}{ x^{n+r} \, y^{n+r} \: dx \, dy}}  \nonumber \\
&=& \sum_{n=0}^\infty{ \left( \int_0^{1^{-}}{x^{n+r} \: dx} \: \int_0^{1^{-}}{ y^{n+r} \: dy}\right)} = \sum_{n=0}^\infty{ \left\{\left[ \frac{x^{n+r+1}}{n+r+1} \right]_0^{1^{-}} \cdot \left[ \frac{y^{n+r+1}}{n+r+1} \right]_0^{1^{-}}\right\}} \nonumber \\
&=& \sum_{n=0}^\infty{ \, \lim_{\varepsilon \rightarrow \, 0^{+}}{ \left( \frac{(1-\varepsilon)^{n+r+1}}{n+r+1} \right)^{\!2}}} = \sum_{n=0}^\infty{ \frac{1}{(n+r+1)^2} } = \sum_{m=1}^\infty{1/(m+r)^2} \, ,
\label{eq:Irr1}
\end{eqnarray}
which readily expands to $\;\sum_{m=1}^{\infty}{1/m^2} \,- \sum_{m=1}^{\,r}{1/m^2}$.
\begin{flushright} $\Box$ \end{flushright}
\end{prova}

For any $\,n \in \mathbb{N}$, let $\,H_n := \sum_{k=1}^n{1/k}\,$ be the $n$-th harmonic number, except for $\,n=0$, for which we define $\,H_0 := 0$.

\begin{lema}[$I_{rs}$]
\label{lem:Irs}
\; Let $r$ and $s$ be non-negative integers, with $\,r \ne s$. Then
\begin{equation*}
\int_0^1\!\!\!\int_0^1{x^r y^s \, \frac{1}{1-x y} ~ dx \, dy} = \frac{H_r -H_s}{r-s} \, .
\label{eq:Irs}
\end{equation*}
\end{lema}

\begin{prova}
\; Let $\,I_{rs} := \int_0^1\!\!\int_0^1{x^r y^s /(1-x y) \: dx \, dy}$.  Again, by expanding the integrand in a geometric series one has
\begin{eqnarray}
I_{rs} &=& \lim_{\varepsilon \rightarrow \, 0^{+}}\int_0^{1-\varepsilon}\!\!\int_0^{1-\varepsilon}{\frac{x^r y^s}{1-x y} \: dx \, dy} \nonumber \\
&=& \int_0^{1^{-}}\!\!\!\int_0^{1^{-}}{x^r y^s \, \sum_{n=0}^\infty{(x y)^n} \: dx \, dy} = \sum_{n=0}^\infty{\int_0^{1^{-}}\!\!\!\int_0^{1^{-}}{ x^{n+r} \, y^{n+s} \: dx \, dy}}  \nonumber \\
&=& \sum_{n=0}^\infty{ \left( \int_0^{1^{-}}{x^{n+r} \: dx} \: \int_0^{1^{-}}{ y^{n+s} \: dy}\right)} = \sum_{n=0}^\infty{ \left\{\left[ \frac{x^{n+r+1}}{n+r+1} \right]_0^{1^{-}} \cdot \left[ \frac{y^{n+s+1}}{n+s+1} \right]_0^{1^{-}}\right\}}  \nonumber \\
&=& \sum_{n=0}^\infty{ \, \lim_{\varepsilon \rightarrow \, 0^{+}}{ \left[\left( \frac{(1-\varepsilon)^{n+r+1}}{n+r+1} \right) \cdot \left( \frac{(1-\varepsilon)^{n+s+1}}{n+s+1} \right)\right]}} \nonumber \\
&=& \sum_{n=0}^\infty{ \frac{1}{n+r+1} \; \frac{1}{n+s+1}} = \sum_{m=1}^\infty{\frac{1}{m+r} \: \frac{1}{m+s}} \, .
\label{eq:Irs1}
\end{eqnarray}
Assume, without loss of generality, that $\,r > s \ge 0$. As the latter series is a telescopic one, we can make
\begin{equation}
\frac{1}{m+r} \: \frac{1}{m+s} = \frac{1}{m+(s+\Delta)} \: \frac{1}{m+s} = \frac{1}{k+\Delta} \: \frac{1}{k} \, ,
\end{equation}
where we substituted $\,k=m+s\,$ and $\,\Delta=r-s$, both being positive integers.  On applying partial fractions decomposition, we easily find
\begin{equation}
 \frac{1}{m+r} \: \frac{1}{m+s} = \frac{1/(r-s)}{k} - \frac{1/(r-s)}{k+\Delta} \, ,
\end{equation}
which leads us to
\begin{equation}
I_{rs} = \frac{1}{r-s} \sum_{k=s+1}^\infty{\!\left( \frac{1}{k} - \frac{1}{k+(r-s)} \right)} = \frac{1}{r-s} \left(\frac{1}{s+1} +\frac{1}{s+2} +\ldots+ \frac{1}{r}\right).
\label{eq:telesc1}
\end{equation}
\begin{flushright} $\Box$ \end{flushright}
\end{prova}

The basic idea for showing that a given real number $\xi$ is irrational is to construct a (infinite) sequence of non-null linear forms (in $\mathbb{Z}$) $\left\{a_n +b_n \, \xi \right\}_{n\ge1}$ which tends to zero as $\,n \rightarrow \infty$. Indeed, if $\,\xi\,$ were rational then the sequence would be bounded away from zero, independently of $\,n$. So, let us build a such sequence.

Let $~d_r\,$ denotes the least common multiple (lcm) of the first $\,r\,$ positive integers, i.e.~$d_r := \mathrm{lcm}{\{1,2,\ldots,r\}}$.

\begin{lema}[$I_{rr}$ as a linear form]
\label{lem:linIrr}
\; For all $\,r \in \mathbb{N}$,
\begin{equation*}
I_{rr} = \zeta(2) -\frac{z_r}{(d_r)^2}
\label{eq:linIrr}
\end{equation*}
for some $\,z_r \in \mathbb{N}^{*}$.  The only exception is $\,r=0$, for which $\,I_{00}=\zeta{(2)}$.
\end{lema}

\begin{prova}
\; For $\,r=0$, we use Lemma~\ref{lem:I0}, which yields $\,I_{00}=\zeta{(2)}$.  For $\,r>0$, from Lemma~\ref{lem:Ir} we know that
\begin{equation}
I_{rr} = \zeta{(2)} -\left(1 +\frac{1}{2^2} + \ldots +\frac{1}{r^2}\right).
\end{equation}
Then, all we need to prove is that
\begin{equation}
\left(d_r\right)^2 \cdot \left(1 +\frac{1}{2^2} + \ldots +\frac{1}{r^2}\right) \in \mathbb{N}^{*}.
\end{equation}
Firstly, note that
\begin{equation}
d_{r^2} \cdot \left(1 +\frac{1}{2^2} + \ldots +\frac{1}{r^2}\right) = d_{r^2} +\frac{d_{r^2}}{2^2} + \ldots +\frac{d_{r^2}}{r^2}
\end{equation}
is a positive integer since $\,d_{r^2} = \mathrm{lcm}\left\{ 1^2,2^2,\ldots,r^2 \right\}\,$ contains all prime factors of the numbers $\,1^2, 2^2, \ldots, r^2$. Secondly, note that $\,d_{r^2} = (d_r)^2$, which is a consequence of the uniqueness of the prime factors decomposition of any positive integer, i.e.~the Fundamental Theorem of Arithmetic, with the constrain of choosing the greater power for each prime factor. In fact, given two positive integers $a$ and $b$, since $\,\mathrm{lcm}\{a,b\} = \prod_{p_i}{p_i^{\max(\alpha_i,\beta_i)}}$, where $\,a=\prod{p_i^{\alpha_i}}\,$ and $\,b=\prod{p_i^{\beta_i}}\,$ are the mentioned prime factors decompositions, with $\alpha_i,\beta_i \ge 0$, then $\,\mathrm{lcm}\{ a^2,b^2 \} = \prod_{p_i}{p_i^{2\,\max(\alpha_i,\beta_i)}} = [\mathrm{lcm}\{a,b\}]^2$. The extension to more than two positive integers is trivial.
\begin{flushright} $\Box$ \end{flushright}
\end{prova}

\begin{lema}[$I_{rs}$ is a positive rational]
\label{lem:linIrs}
\; For all $~r,s \in \mathbb{N}$, $r \ne s$,
\begin{equation*}
I_{rs} = \frac{z_{rs}}{(d_r)^2}
\label{eq:linIrs}
\end{equation*}
for some $\:z_{rs} \in {\mathbb{N}}^{*}$.
\end{lema}

\begin{prova}
\; Let us assume, without loss of generality, that $\,r > s \ge 0$.  From Lemma~\ref{lem:Irs}, we know that
\begin{equation}
I_{rs} = \frac{1}{r-s} \, \left(\frac{1}{s+1} +\frac{1}{s+2} +\ldots+ \frac{1}{r}\right),
\end{equation}
which means that
\begin{eqnarray}
I_{rs} \cdot (d_r)^2 = \frac{(d_r)^2}{r-s} \cdot \left(\frac{1}{s+1} +\frac{1}{s+2} +\ldots+ \frac{1}{r}\right)  \nonumber \\
= \frac{d_r}{r-s} \, \left(\frac{d_r}{s+1} +\frac{d_r}{s+2} +\ldots+ \frac{d_r}{r}\right).
\end{eqnarray}
Clearly, the last expression is the product of two positive integers since $\,d_r\,$ is a multiple of $\,r-s$, which is a positive integer smaller than (or equal to) $\,r$, as well as a multiple of every integer from $\,s+1\,$ to $\,r$. Therefore, $\,I_{rs} \cdot (d_r)^2 = z_{rs}\,$ is a positive integer.
\begin{flushright} $\Box$ \end{flushright}
\end{prova}

Summarizing, given $\,r,s \in \mathbb{N}$ one has $\,I_{rs} \in \delta_{r s} \, \zeta{(2)} \pm \, \mathbb{N}/d_r^{\,2}$, where $\,\delta_{ij}\,$ is the Kronecker delta and the minus sign is for $\,r = s$.  Therefore, for any polynomials with integer coefficients $\,R_n(x) = \sum_{k=0}^n{a_k\,x^k}\,$ and $\,S_n(y) = \sum_{k=0}^n{b_k\,y^k}$, one has
\begin{eqnarray}
\int_0^1\!\!\!\int_0^1{R_n(x) \: S_n(y) \: \frac{1}{1-x y} \: dx \, dy} = \int_0^1\!\!\!\int_0^1{\sum_{r=0}^n{a_r\,x^r} \cdot \sum_{s=0}^n{b_s\,y^s} \: \frac{1}{1-x y} \: dx \, dy}  \nonumber \\
= \sum_{r=0}^n{\sum_{s=0}^n{ a_r \, b_s \, \int_0^1\!\!\!\int_0^1{\,x^r \, y^s \: \frac{1}{1-x y} \: dx \, dy}}} = \sum_{r=0}^n{\sum_{s=0}^n{ a_r \, b_s \, I_{rs}}} \, , \quad
\label{eq:previoussec}
\end{eqnarray}
which clearly belongs to $\; \mathbb{Z} \: \zeta{(2)} + \mathbb{Z}/d_n^{\:2}$.

\subsection{Legendre-type polynomials and the prime number theorem}

Let us take into account, in the Beukers integrals of the previous section, the following Legendre-type polynomials (normalized in the interval $[0,1]$):\footnote{Note that $\,P_n(0)=1\,$ and $\,P_n(1)=(-1)^n$. In fact, this is a particular case of a general symmetry rule, namely $\,P_n(1-x) = (-1)^n \, P_n(x)$.}
\begin{equation}
P_{n}(x) := \frac{1}{n!} \: \frac{d^n}{dx^n}\left[x^n\,(1-x)^n \right].
\label{eq:defPn}
\end{equation}
As pointed out in Ref.~\cite{Beukers80}, these polynomials can be written in the equivalent form\footnote{For example: $\,P_0(x)=1$, $\,P_1(x)=1-2x$, $\,P_2(x)=1-6x+6x^2$, etc.}
\begin{equation}
P_{n}(x) = \sum_{k=0}^n{(-1)^k \, \binom{n}{k} \, \binom{n+k}{n} \, x^k} \, ,
\end{equation}
which makes it clear that $\,P_n(x)\,$ has integer coefficients for all $\,n \in \mathbb{N}$,  so the linear forms in $\mathbb{Z}$ obtained in the end of the previous section do apply when we put $\,R_n(x) = P_n(x)\,$ in Eq.~\eqref{eq:previoussec}.

The choice of Legendre-type polynomials comes from the possibility of performing integration by parts easily, as describes the following lemma.

\begin{lema}[Integration by parts with $\,P_n(x)\,$]
\label{lem:intPnx}
\; For all $\,n \in \mathbb{N}\,$ and $\,f: [0,1] \rightarrow \mathbb{R}\,$ of class $\,\mathcal{C}^n$, one has
\begin{equation*}
\int_0^1{P_n(x) \, f(x) \: dx} = \frac{(-1)^n}{n!} \: \int_0^1{x^n \, (1-x)^n ~ \frac{d^n f}{dx^n} ~ dx} \, .
\label{eq:intIn}
\end{equation*}
\end{lema}

\begin{prova}
\; The proof is a sequence of integration by parts, but it suffices to make the first one. Given $\,n \in \mathbb{N}\,$ and $\,f: [0,1] \rightarrow \mathbb{R}\,$ of class $\,\mathcal{C}^n$, from the definition of $\,P_n(x)\,$ in Eq.~\eqref{eq:defPn}, it follows that
\begin{eqnarray}
I_n := \int_0^1{P_n(x) \, f(x) \: dx} = \int_0^1{\frac{1}{n!} \: \frac{d^n}{dx^n}\left[x^n\,(1-x)^n \right] \: f(x) ~ dx}  \nonumber \\
= \frac{1}{n!} \, \int_0^1{\frac{d}{dx} \left\{ \frac{d^{n-1}}{dx^{n-1}}\left[x^n\,(1-x)^n \right] \right\} \, f(x) ~ dx} \, .
\end{eqnarray}
Since $\,\int{u \, dv} = u \, v -\int{v \, du}$, then let us choose $\,u=f(x)\,$ and $\,dv = d/dx\{\ldots\} \, dx$. With this choice, $du = f'(x) \, dx\,$ and $\,v = \{\ldots\}$, so
\begin{eqnarray}
n! \: I_n = \left[f(x) \, \frac{d^{n-1}}{dx^{n-1}}\left(x^n \, (1-x)^n \right)\right]_0^1 -\int_0^1{\frac{d^{n-1}}{dx^{n-1}}\left[x^n\,(1-x)^n\right] \: f'(x)} \: dx  \nonumber \\
= \left[f(x) \, \sum_{k=0}^{n-1}{\binom{n-1}{k}\,(x^n)^{(k)} \, \left[(1-x)^n\right]^{(n-1-k)} } \right]_0^1 -\int_0^1{\left[x^n\,(1-x)^n\right]^{(n-1)} \, f'(x)} \: dx \, , \quad
\end{eqnarray}
where we have made use of the generalized Leibnitz rule for the higher derivatives of the product of two functions, i.e. $\,(f \cdot g)^m = \sum_{k=0}^m{\binom{m}{k}\, f^{(k)} \, g^{(m-k)}}$.  On calculating these derivatives, one finds
\begin{eqnarray}
n! \: I_n = \left[f(x) \, \sum_{k=0}^{n-1}{\frac{(n-1)!}{k!\,(n-1-k)!} \, \frac{n!}{(n-k)!} \, x^{n-k} \, \frac{n!}{(k+1)!} \, (-1)^k \, (1-x)^{k+1} } \right]_0^1  \nonumber \\
-\int_0^1{\left[x^n\,(1-x)^n\right]^{(n-1)} \, f'(x)} \: dx  \nonumber \\
= \, \left[ f(1) \times 0 - f(0) \times 0 \right] -\int_0^1{\left[x^n\,(1-x)^n\right]^{(n-1)} \, f'(x)} \: dx  \nonumber \\
= -\int_0^1{\left[x^n\,(1-x)^n\right]^{(n-1)} \, f'(x)} \: dx \, , \;
\end{eqnarray}
where the zeros in $\,f(1) \times 0\,$ and $\,f(0) \times 0\,$ come from the presence of factors $\,(1-x)\,$ and $\,x$, respectively, in every terms of the finite sum on $k$. Then, the overall result of the first integration by parts is the transference of a derivative $\,d/dx\,$ from $\,P_n(x)\,$ to the function $f(x)$ and a change of sign. Of course, each further integration by parts will produce the same effect, so it is easy to deduce that $\, n! ~ I_n = (-1)^n \int_0^1{x^n\,(1-x)^n \: f^{(n)}(x) \; dx}$.
\begin{flushright} $\Box$ \end{flushright}
\end{prova}

The appearance of $\,d_n\,$ in the linear forms we are treating here will demand, during the proof of the main theorem, a `good' upper bound for it.

\begin{lema}[Upper bound for $d_n$]
\label{lem:dn}
\; Let $\,n\,$ be a positive integer. Define $\,\pi{(n)}\,$ as the number of primes less than (or equal to) $\,n$. Then, $\,d_n \le n^{\pi{(n)}}\,$ and $\,n^{\pi{(n)}} \sim e^n$.
\end{lema}

\begin{prova}
\; For all $\,n \in \mathbb{N}^{*}$, $d_n = \mathrm{lcm}\{1,2,\ldots,n\}\,$ is formed by multiplying together all primes $\,p \le n$ with the greatest possible exponents $\,m\,$ such that $\,p^m \le n$. Therefore
\begin{equation}
d_n = \prod_{p \le n}{p^m} \: ,
\end{equation}
where $\,m := \max_{k \in \mathbb{N}}\{p^k \le n\}$. Of course, the largest exponent $\,m\,$ is so that $\,p^m \le n$. On taking the logarithm (with basis $p$) on both sides of this inequality, one finds $\,m \le \log_p{n}$. Since $\,m \in \mathbb{N}\,$ is to be maximal, then we must take $\,m = \lfloor \log_p{n} \rfloor = \lfloor \ln{n}/\ln{p} \rfloor$. This implies that
\begin{equation}
d_n = \prod_{p \le n}{p^{\lfloor \ln{n}/\ln{p} \rfloor}} \: .
\end{equation}
In fact, $\,m = \lfloor \ln{n}/\ln{p} \rfloor \le \ln{n}/\ln{p}\:$ implies that
\begin{equation}
p^{\left\lfloor \frac{\ln{n}}{\ln{p}} \right\rfloor} \le p^{\frac{\ln{n}}{\ln{p}}} = \left(e^{\ln{p}}\right)^{\frac{\ln{n}}{\ln{p}}} = e^{\ln{n}} = n \, ,
\end{equation}
and then
\begin{equation}
d_n = \prod_{p \le n}{p^{\lfloor \ln{n}/\ln{p} \rfloor}} \le \prod_{p \le n}{n} = n^{\pi{(n)}}\, .
\end{equation}

From the prime number theorem (PNT), we know that the asymptotic behavior of the function $\,\pi{(n)}\,$ is the same of the function $\,n/\ln{n}$, for sufficiently large values of $n$. Then,
\begin{equation}
\pi{(n)} \sim \frac{n}{\ln{n}} \: \Longrightarrow \: n^{\pi{(n)}} \sim n^{n/\ln{n}} = \left( e^{\ln{n}} \right)^{n/\ln{n}} = e^n \, ,
\end{equation}
so $\: n^{\pi{(n)}} \sim e^n$.
\begin{flushright} $\Box$ \end{flushright}
\end{prova}

Now we have all ingredients in hands to prove the first irrationality result.

\subsection{First main result}

\begin{teo}[$\zeta{(2)} \not\in \mathbb{Q}$]
\label{teo:main}
\quad The number $\,\zeta{(2)}\,$ is irrational.
\end{teo}

\begin{prova}
\quad  In Lemma~\ref{lem:intPnx}, choose $\,f(x) = \int_0^1{(1-y)^n/(1-xy) ~ dy}$.  From Lemmas~\ref{lem:linIrr} and~\ref{lem:linIrs}, for all $\,n \in \mathbb{N}$ we have
\begin{equation}
I_n = \int_0^1{P_n(x) \, f(x) \: dx} = \int_0^1\!\!\!\int_0^1{P_n(x) \, \frac{(1-y)^n}{1 -x y} \: dy \, dx} \, \in \, \mathbb{Z} \, \zeta{(2)} + \frac{\mathbb{Z}}{d_n^{\:2}} \, .
\end{equation}
In particular, $I_0 = \int_0^1{P_0(x) \, f(x) \: dx} = \int_0^1{1 \, \left[\int_0^1{\frac{1}{1 -x y}~dy} \right] \, dx} = \int_0^1\!\!\int_0^1{\frac{1}{1 -x y}~dx\,dy} = \zeta(2)$, as shown in Lemma~\ref{lem:I0}.  Hence,
\begin{equation}
I_n = \frac{a_n}{d_n^{\:2}} + b_n\,\zeta{(2)} \, ,
\label{eq:Inanbn}
\end{equation}
for some integers $\,a_n$ and $\,b_n$.

On the other hand,
\begin{equation}
|I_n| = \left| \int_0^1{P_n(x) \cdot \left(\int_0^1{\frac{(1-y)^n}{1-xy} \: dy} \right) dx} \right| ,
\end{equation}
which, from Lemma~\ref{lem:intPnx}, becomes
\begin{eqnarray}
|I_n| = \left| \frac{(-1)^n}{n!} \, \int_0^1{x^n \, (1-x)^n \, \frac{d^n}{d x^n} \left(\int_0^1{\frac{(1-y)^n}{1-xy} \: dy} \right) dx} \right|  \nonumber \\
= \frac{1}{n!} \, \left| \int_0^1{x^n \, (1-x)^n \, \int_0^1{\frac{\partial^n}{\partial x^n}\left(\frac{(1-y)^n}{1-xy} \right) \: dy} \: dx} \right|  \nonumber \\
= \frac{1}{n!} \, \left| \int_0^1{x^n \, (1-x)^n \, \int_0^1{(1-y)^n \frac{\partial^n}{\partial x^n}\left(1-y x\right)^{-1} \: dy} \: dx} \right|  \nonumber \\
= \frac{1}{n!} \, \left| \int_0^1{x^n \, (1-x)^n \, \int_0^1{(1-y)^n \, \frac{n! \: y^n}{(1-y x)^{n+1}} \: dy} \: dx} \right|  \nonumber \\
= \int_0^1\!\!\!\int_0^1{\frac{x^n \, (1-x)^n \, y^n \, (1-y)^n}{(1- x y)^{n+1}} ~ dx \: dy} \, .
\end{eqnarray}
Note that the integrand in the latter integral is positive over all points of $\,(0,1)^2$, being null only at the boundaries of $[0,1]^2$, except at the point $(1,1)$, where it is an indeterminate form of the kind `$0/0$', though it tends to zero as $(x,y) \rightarrow (1^{-},1^{-})$ for all $n>1$.\footnote{For $n=1$, this integrand tends to $1/4$ as $(x,y) \rightarrow (1^{-},1^{-})$, as the reader can easily check. For $n=0$, the integral reduces to $I_{00} = \zeta(2)$, as seen in Lemma~\ref{lem:I0}, which of course is positive. }  Therefore, $\,|I_n| > 0$, $\forall \, n \in \mathbb{N}$.\footnote{In fact, Lemma~1 of Ref.~\cite{Alladi1979} shows that, given a continuous function $f(x)$, if $\,\int_0^1{P_n(x) \: f(x)~dx} = 0\,$ for all $n \ge N$, then $f(x)$ necessarily is a polynomial of degree $N$ or smaller.}

On searching for a suitable upper bound for $\,|I_n|$, all we need to do is
\begin{eqnarray}
|I_n| = \int_0^1\!\!\!\int_0^1{\left[\frac{\,x \, (1-x) \, y \, (1-y)}{1- x y}\right]^n  \frac{1}{1 -x y} ~ dx \: dy}  \nonumber \\
\le \int_0^1\!\!\!\int_0^1{ \max_{[0,1)^2}\left\{\left[\frac{\,x \, (1-x) \, y \, (1-y)}{1- x y}\right]^n \right\}  \frac{1}{1 -x y} ~ dx \: dy}  \nonumber \\
= \left\{\max_{[0,1)^2}\left[\frac{\,x \, (1-x) \, y \, (1-y)}{1- x y}\right]\right\}^n  \, \int_0^1\!\!\!\int_0^1{\frac{1}{1 -x y} ~ dx \: dy}  \nonumber \\
= \left\{\max_{[0,1)^2}\left[\frac{\,x \, (1-x) \, y \, (1-y)}{1- x y}\right]\right\}^n  \, \zeta{(2)}\, .
\end{eqnarray}
The above maximum can be determined analytically. For this, let
\begin{equation*}
g(x,y) := \frac{x \, (1-x) \, y \, (1-y)}{1- x y}
\end{equation*}
be the function we have to maximize over $[0,1)^2$. Firstly, note that $\,g(x,y) > 0 \,$ throughout the open $(0,1)^2\,$ and $\,g(x,y)=0\,$ in all points of the boundary of $[0,1)^{2\,}$.\,\footnote{Note that $\,g(x,y) \rightarrow 0\,$ as $\,(x,y) \rightarrow (1^{-},1^{-})$, so we do not need to worry about this point.}  Now, since both $\,\partial^2 g/\partial x^2 = -2 y (1-y)^2/(1-xy)^3 \,$ and $\,\partial^2 g/\partial y^2 = -2 x (1-x)^2/(1-xy)^3 \,$ are negative for all $\,(x,y) \in (0,1)^2\,$ and $\,g(x,y) = g(y,x)$, then the maximum is unique and has the form $\,g(t,t)$, for some $\,t \in (0,1)$. Since $\,g(t,t) = (t -t^2)^2/(1-t^2)$, then $\,d g / d t = 0\,$ leads to $\,t^2+t-1=0$, which has two real solutions, namely $\,t = (-1 \pm \sqrt{5})/2$. The only solution in the open $(0,1)$ is $\,\Phi = (\sqrt{5}-1)/2$, which is the inverse of the golden ratio $\,\phi = (\sqrt{5}+1)/2$. This yields $\,\max_{[0,1)^2}{[g(x,y)]} = g(\Phi,\Phi) = \Phi^5$, thus $\,|I_n| \le \Phi^{5n} \, \zeta{(2)}\,$ for all $\,n \in \mathbb{N}$. Since $\,\Phi < 1$, this shows that $\,|I_n| \rightarrow 0\,$ as $n \rightarrow \infty$.\footnote{As pointed out at the end of Sec.~3 of Ref.~\cite{Alladi1980}, $|I_n| \rightarrow 0\,$ as $n \rightarrow \infty\,$ for any function $f(x) \in \mathcal{L}^2_{[0,1]}$, as follows from a well-known property of Legendre polynomials.} From Eq.~\eqref{eq:Inanbn}, one has
\begin{equation}
0 < \left| \frac{a_n}{d_n^{\:2}} + b_n \, \zeta{(2)} \right| \le \Phi^{5n} \: \zeta{(2)} \, ,
\end{equation}
which means that
\begin{equation}
0 < \left| a_n + b_n \, d_n^{\:2} \: \zeta{(2)} \right| \le {d_n^{\:2}} \,~ \Phi^{5n} \: \zeta{(2)} \, .
\end{equation}
From Lemma~\ref{lem:dn}, we know that $\,d_n^{\:2} \le (n^{\pi{(n)}})^2\,$ and $\, (n^{\pi{(n)}})^2 \sim \left(e^n\right)^2 = (e^2)^n$. Since $\,e^2 < 8$, then, for sufficiently large values of $n$, $d_n^{\:2} < 8^n$, which leads us to
\begin{equation}
0 < \left| a_n + c_n \: \zeta{(2)} \right| < 8^n ~ \left(\Phi^5\right)^n \; \zeta{(2)} = \left(8\, \Phi^5 \right)^n \: \zeta{(2)} \, ,
\end{equation}
where $\: c_n = b_n \, d_n^{\:2} \,$ is an integer.  Now, assume, towards a contradiction, that $\,\zeta{(2)}  \in \mathbb{Q}$. Since $\zeta{(2)}>0$, then $\,\zeta{(2)}\,$ can be written in the form $\,p/q$, with $p$ and $q$ being coprime positive integers.  From the above inequality, one has
\begin{equation}
0 < \left| a_n + c_n \: \frac{p}{q} \right| < \left(8\, \Phi^5 \right)^n \: \frac{p}{q} \, .
\end{equation}
On searching for a contradiction, note that $\: 8 \, \Phi^5 = 0.7213\ldots < 0.75 = 3/4$, so
\begin{equation}
0 < \left| q \, a_n + p \: c_n \right| < p \left( \frac{3}{4} \right)^{\!n} .
\end{equation}
Since $\,\left| q \, a_n + p \: c_n \right|\,$ is a positive integer, then $\,\left| q \, a_n + p \: c_n \right| \ge 1\,$ for all $\,n \in \mathbb{N}$. However, the fact that $\,\left| q \, a_n + p \: c_n \right| < p \, \left( 3/4 \right)^n\,$ forces $\, \left| q \, a_n + p \: c_n \right| \,$ to be less than $\,1\,$ for sufficiently large values of $\,n$, and we have a contradiction.\footnote{More precisely, $\,\left| q \, a_n + p \: c_n \right| < 1\:$ for all $\:n > \ln{p\,}/\ln{(4/3)}$.}  Therefore, $\,\zeta{(2)}\:$ cannot be a positive rational number.
\begin{flushright} $\Box$ \end{flushright}
\end{prova}

\section{Irrationality of $\zeta{(3)}$}

\subsection{Preliminaries}

As for $\zeta{(2)}$, we begin with some lemmas on unit square integrals.

\begin{lema}[An unit square integral for $\,\zeta{(3)}\,$]
\label{lem:I0z3}
\begin{equation*}
- \int_0^1\!\!\!\int_0^1{\frac{\ln{(x y)}}{1-x y} \: dx \, dy} = 2 \, \zeta{(3)} \, .
\label{eq:I0z3}
\end{equation*}
\end{lema}

\begin{prova}
\; Define $\,J_{00} := - \int_0^1\!\!\int_0^1{\ln{(x y)}/(1-x y) \: dx \, dy}$.  It follows that
\begin{eqnarray}
J_{00} &=& - \lim_{\varepsilon \rightarrow \, 0^{+}}\int_0^{1-\varepsilon}\!\!\int_0^{1-\varepsilon}{\frac{\ln{(x y)}}{1-x y} ~ dx \, dy} = - \int_0^{1^{-}}\!\!\!\int_0^{1^{-}}{\sum_{n=0}^\infty{(x y)^n \, \ln{(x y)}} \: dx \, dy} \nonumber \\
&=& - \sum_{n=0}^\infty{\int_0^{1^{-}}\!\!\!\int_0^{1^{-}}{\!x^n \, y^n \, \left(\ln{x} +\ln{y}\right)} \: dx \, dy} = - 2 \sum_{n=0}^\infty{\left( \int_0^{1^{-}}{\!x^n \, \ln{x} \: dx} \, \cdot \, \int_0^{1^{-}}{\!y^n \: dy} \right) } \, . \; \qquad
\end{eqnarray}
Integration by parts yields, apart from an arbitrary constant of integration, $\int{x^n\,\ln{x} \, dx} = x^{n+1}\,\ln{x}/(n+1) -x^{n+1}/(n+1)^2$, so
\begin{eqnarray}
J_{00} &=& - 2 \sum_{n=0}^\infty{\left[\frac{x^{n+1}\,\ln{x}}{n+1} -\frac{x^{n+1}}{(n+1)^2}\right]_0^{1^{-}} \times \left[\frac{y^{n+1}}{n+1}\right]_0^{1^{-}}}  \nonumber \\
&=& - 2 \sum_{n=0}^\infty{ \, \lim_{\varepsilon \rightarrow \, 0^{+}}{ \left[ \frac{(1-\varepsilon)^{n+1} \, \ln{(1-\varepsilon)}}{n+1} -\frac{(1-\varepsilon)^{n+1}}{(n+1)^2} -\frac{\varepsilon^{n+1} \, \ln{\varepsilon}}{n+1} +0 \right] \times \left[ \frac{(1-\varepsilon)^{n+1}}{n+1} -0 \right]}}  \nonumber \\
&=& - 2 \sum_{n=0}^\infty{ \left( \frac{\ln{1}}{n+1} -\frac{1}{(n+1)^2} -0 \right) \times \frac{1}{n+1} } = 2 \sum_{n=0}^\infty{ \frac{1}{(n+1)^3} } = 2 \sum_{m=1}^\infty{ \frac{1}{m^3}} \, .
\end{eqnarray}
\begin{flushright} $\Box$ \end{flushright}
\end{prova}

\begin{lema}[$J_{rr}$]
\label{lem:Irz3}
\; For all integers $\,r>0$
\begin{equation*}
- \int_0^1\!\!\!\int_0^1{x^r y^r \: \frac{\,\ln{(x y)}}{1-x y} ~ dx \, dy} = 2 \, \zeta{(3)} - 2 \sum_{m=1}^r{\frac{1}{m^3}} \, .
\label{eq:Irz3}
\end{equation*}
\end{lema}

\begin{prova}
\; Define $\:J_{rr} := - \int_0^1\!\!\int_0^1{x^r y^r \, \ln{(x y)}/(1-x y) \: dx\,dy}$. On substituting $\,r\,$ by $\,r +t\,$ in Eq.~\eqref{eq:Irr1}, $\,t\,$ being a positive real, one finds
\begin{equation}
\int_0^1\!\!\!\int_0^1{x^{r+t} \, y^{r+t}\, \frac{1}{1-x y} ~ dx\,dy} = \sum_{m=1}^\infty{\frac{1}{(m+r+t)^2}} \, .
\end{equation}
On taking the derivative with respect to $\,t\,$ on both sides, one has\footnote{Again, the interchange of limits, sums, derivatives and integrals is fully justified in the proof of Theorem 2.1 of Ref.~\cite{Hadji2002}.}
\begin{eqnarray}
& & \frac{d}{d t} \int_0^1\!\!\!\int_0^1{x^{r+t} y^{r+t}\, \frac{1}{1-x y} ~dx\,dy} = \frac{d}{d t} \sum_{m=1}^\infty{\frac{1}{(m+r+t)^2}}  \nonumber \\
& \Longrightarrow & \, \int_0^1\!\!\!\int_0^1{\frac{1}{1-x y} ~ \frac{\partial}{\partial t} \left(x \, y\right)^{t+r} \: dx\,dy} = \sum_{m=1}^\infty{\frac{d}{d t} \left[\frac{1}{(t +r +m)^2}\right]}  \nonumber \\
& \Longrightarrow & \, \int_0^1\!\!\!\int_0^1{\frac{\,(x \, y)^r}{1-x y} \: \left(x \, y\right)^t \, \ln{(x y)}~ dx\,dy} = - 2 \sum_{m=1}^\infty{\frac{1}{(t +r +m)^3}} \, .
\end{eqnarray}
On putting $\,t=0$, this reduces to
\begin{equation}
\int_0^1\!\!\!\int_0^1{\frac{(x \, y)^r}{1-x y} \: \ln{(x y)}~ dx\,dy} = -J_{rr} = - \,2 \sum_{m=1}^\infty{\frac{1}{(r +m)^3}} \, .
\end{equation}
The last sum readily expands to $ ~ \sum_{m=1}^{\infty}{1/m^3} \,- \sum_{m=1}^{\,r}{1/m^3}$, which completes the proof.
\begin{flushright} $\Box$ \end{flushright}
\end{prova}

In view to extend the definition of $\,H_n\,$ to denominators with exponent $2$, define $\:H_n^{(2)} := \sum_{k=1}^n{1/k^2}$, except for $\,n=0$, for which we define $\:H_0^{(2)}:=0$.

\begin{lema}[$J_{rs}$]
\label{lem:Irsz3}
\; Let $r$ and $s$ be non-negative integers, with $\,r \ne s$. Then
\begin{equation*}
- \int_0^1\!\!\!\int_0^1{x^r \, y^s \: \frac{\,\ln{(x y)}}{1-x y} ~ dx \, dy} = \frac{H_r^{(2)} -H_s^{(2)}}{r-s} \, .
\label{eq:Irsz3}
\end{equation*}
\end{lema}

\begin{prova}
\; Define $\:J_{rs} := - \int_0^1\!\!\int_0^1{x^r \, y^s \, \ln{(x y)}/(1-x y) \: dx\,dy}$.  From the last line of Eq.~\eqref{eq:Irs1}, one has
\begin{equation}
\int_0^1\!\!\!\int_0^1{x^{r+t} \, y^{s+t}\, \frac{1}{1-x y} ~ dx\,dy} = \sum_{m=1}^\infty{\frac{1}{(m+r+t)\,(m+s+t)}} \, .
\end{equation}
Again, on taking the derivative with respect to $\,t\,$ on both sides, one finds
\begin{equation}
\int_0^1\!\!\!\int_0^1{\frac{x^r\,y^s}{1-x y} ~ \frac{\partial}{\partial t} \left(x \, y\right)^t \: dx\,dy} = \sum_{m=1}^\infty{\frac{d}{d t} \left[\frac{1}{(m+r+t)\,(m+s+t)}\right]} \, . \quad
\end{equation}
On assuming, without loss of generality, that $\,r>s$, one finds
\begin{eqnarray}
\int_0^1\!\!\!\int_0^1{\frac{\,x^r \, y^s}{1-x y} \: \left(x \, y\right)^t \, \ln{(x y)}~ dx\,dy} = \sum_{m=1}^\infty{\frac{d}{d t} \left[\frac{1}{r-s} \left(\frac{1}{m+s+t} -\frac{1}{m+r+t} \right)\right]} \nonumber \\
= \sum_{m=1}^\infty{\frac{1}{r-s} ~ \frac{d}{d t} \left(\frac{1}{m+s+t} -\frac{1}{m+r+t} \right)} = \sum_{m=1}^\infty{\frac{1}{r-s} \left[\frac{-1}{(m+s+t)^2} -\frac{(-1)}{(m+r+t)^2} \right]} \nonumber \\
= \frac{-1}{r-s} \, \sum_{m=1}^\infty{\left[\frac{1}{(m+s+t)^2} -\frac{1}{(m+r+t)^2} \right]} \, . \quad
\end{eqnarray}
For $\,t=0$, this reduces to
\begin{equation}
\int_0^1\!\!\!\int_0^1{\frac{\,x^r \, y^s}{1-x y} \: \ln{(x y)}~ dx\,dy} = -J_{rs} = \frac{-1}{r-s} \, \sum_{m=1}^\infty{\left[\frac{1}{(m+s)^2} -\frac{1}{(m+r)^2} \right]} \, .
\end{equation}
The last sum telescopes to $\; 1/(s+1)^2 +\ldots +1/r^2 = H_r^{(2)} -H_s^{(2)}$.
\begin{flushright} $\Box$ \end{flushright}
\end{prova}

Now, let us make use of the above integrals to construct linear forms in $\mathbb{Z}$ involving $\,\zeta{(3)}$.

\begin{lema}[$J_{rr}$ as a linear form]
\label{lem:linIrrz3}
\; For all $\,r \in \mathbb{N}$,
\begin{equation*}
J_{rr} = 2\,\zeta{(3)} -\frac{z_r}{(d_r)^3}
\label{eq:linIrrz3}
\end{equation*}
for some $\,z_r \in \mathbb{N}^{*}$.  The only exception is $\,r=0$, for which $\,J_{00} = 2 \, \zeta{(3)}$.
\end{lema}

\begin{prova}
\; For $\,r=0$, we use Lemma~\ref{lem:I0z3}, which yields $\,J_{00} = 2 \, \zeta{(3)}$.  For $\,r>0$, from Lemma~\ref{lem:Irz3} we know that
\begin{equation}
J_{rr} = 2\,\zeta{(3)} - 2 \left(1 +\frac{1}{2^3} + \ldots +\frac{1}{r^3}\right).
\end{equation}
Then, all we need to prove is that
\begin{equation}
\left(d_r\right)^3 \cdot \left(1 +\frac{1}{2^3} + \ldots +\frac{1}{r^3}\right) \in \mathbb{N}^{*}.
\end{equation}
Note that
\begin{equation}
d_{r^3} \cdot \left(1 +\frac{1}{2^3} + \ldots +\frac{1}{r^3}\right) = d_{r^3} +\frac{d_{r^3}}{2^3} + \ldots +\frac{d_{r^3}}{r^3}
\end{equation}
is a positive integer since $\,d_{r^3} = \mathrm{lcm}\left\{ 1^3, 2^3, \ldots, r^3 \right\}\,$ contains all prime factors of the numbers $\,1^3, 2^3, \ldots, r^3$. As in the proof of Lemma~\ref{lem:linIrr}, it is easy to see that $\,d_{r^3} = (d_r)^3$, which completes the proof.
\begin{flushright} $\Box$ \end{flushright}
\end{prova}

\begin{lema}[$J_{rs}$ is a positive rational]
\label{lem:linIrsz3}
\; For all $~r,s \in \mathbb{N}$, $r \ne s$,
\begin{equation*}
J_{rs} = \frac{z_{rs}}{(d_r)^3}
\label{eq:linIrsz3}
\end{equation*}
for some $\:z_{rs} \in {\mathbb{N}}^{*}$.
\end{lema}

\begin{prova}
\; Let us assume, without loss of generality, that $\,r > s \ge 0$.  From Lemma~\ref{lem:Irsz3}, we know that
\begin{equation}
J_{rs} = \frac{H_r^{(2)} -H_s^{(2)}}{r-s} = \frac{1}{r-s} \, \left[\frac{1}{\,(s+1)^2} +\frac{1}{(s+2)^2} +\ldots+ \frac{1}{r^2} \, \right] ,
\end{equation}
which means that
\begin{eqnarray}
J_{rs} \cdot (d_r)^3 = \frac{(d_r)^3}{r-s} \, \left(\frac{1}{(s+1)^2} +\frac{1}{(s+2)^2} +\ldots +\frac{1}{r^2}\right)  \nonumber \\
= \frac{d_r}{r-s} \, \left( \frac{(d_r)^2}{s+1} +\frac{(d_r)^2}{s+2} +\ldots +\frac{(d_r)^2}{r} \right).
\end{eqnarray}
Clearly, the last expression is the product of two positive integers since $\,d_r\,$ is a multiple of $\,r-s$, which is a positive integer smaller than (or equal to) $\,r$, and $\,(d_r)^2 = d_{r^2}\,$ is a multiple of each element of $\,\{(s+1)^2, (s+2)^2, \ldots, r^2\}$. Therefore, $\,J_{rs} \cdot (d_r)^3 = z_{rs}\,$ is a positive integer.
\begin{flushright} $\Box$ \end{flushright}
\end{prova}

Summarizing, $\,J_{rs} \in 2 \, \delta_{r s} \, \zeta{(3)} \pm \, \mathbb{N}/d_r^{\:3}$, for all $\,r,s \in \mathbb{N}$, the minus sign being for $\,r = s$.  Analogously to what we have done for $\zeta{(2)}$, given any two polynomials with integer coefficients $\,R_n(x)\,$ and $\,S_n(y)$, one has
\begin{equation}
\int_0^1\!\!\!\int_0^1{R_n(x) \: S_n(y) \: \frac{\ln{(x y)}}{\,1-x y\,} ~ dx \, dy} \: \in \: \mathbb{Z} \: \zeta{(3)} + \mathbb{Z}/d_n^{\:3} \, .
\label{eq:prevZ3}
\end{equation}
This linear form in $\mathbb{Z}$ of course holds for $P_n(x)$ and $P_n(y)$, as they have only integer coefficients.

We need two more results on integrals for proofing our second main result.

\begin{lema}[An useful substitution]
\label{lem:intLn}
\begin{equation*}
\int_0^1{\frac{1}{\,1-(1 -v)\,z\,} ~dz} = - \, \frac{\;\ln{v}\:}{\:1 -v\,} \: .
\label{eq:intLn}
\end{equation*}
\end{lema}

\begin{prova}
Substitute $\,(1 -v)\,z = u\,$ in the integral, with $\, d u = (1 -v)\,dz$. This yields
\begin{eqnarray}
\int_0^1{\frac{1}{\,1-(1 -v)\,z\,} ~dz} = \int_0^{1 -v}{\frac{1}{\,(1 -u) \, (1 -v)\,} ~du} = \frac{1}{1 -v} \, \int_0^{1-v}{\frac{1}{\,1 -u\,} ~du} \nonumber \\
= \frac{1}{1 -v} \, \left[- \ln{(1 -u)}\right]_0^{1 -v} = - \frac{1}{1 -v} \left[ \, \ln{v} -\ln{1} \, \right] = - \frac{\ln{v}}{\,1 -v\,}  \, . \quad
\end{eqnarray}
\begin{flushright} $\Box$ \end{flushright}
\end{prova}

\begin{lema}[A partial fraction integration]
\label{lem:intParcial}
\; Given $s,t \in \mathbb{R}_{(0,1)}$, the following equality holds:
\begin{equation*}
\int_0^1{\frac{1}{\,1-[1 -(1-s)\,t]\,u\,} ~du} = \int_0^1{\frac{1}{\,[1-(1-u)\,s] \: [1 -(1-t)\,u]\,} ~du} \, .
\label{eq:intParcial}
\end{equation*}
\end{lema}

\begin{prova}
On putting $\:v = (1-s)\,t\:$ in Lemma~\ref{lem:intLn}, one finds
\begin{equation}
\int_0^1{\frac{1}{\,1-(1 -(1-s)t)\,z\,} ~dz} = - \,\frac{\;\ln{[(1-s)\,t]}\:}{\:1 -(1-s)\,t\,} \, .
\end{equation}
The other integral can be solved by making use of the following partial fraction decomposition:
\begin{equation*}
\frac{1}{\,[1-(1-u)\,s] \: [1 -(1-t)\,u]} = \frac{1}{1-(1-s)\,t} \left[ \frac{s}{1-(1-u)\,s} -\frac{1-t}{1-(1-t)\,u} \right] ,
\end{equation*}
which implies that
\begin{eqnarray}
\int_0^1{\frac{1}{\,[1-(1-u)\,s] \: [1 -(1-t)\,u]\,} ~du} \nonumber \\
= \frac{1}{1-(1-s)\,t} \, \int_0^1{\left[ \frac{s}{1-(1-u)\,s} -\frac{1-t}{1-(1-t)\,u} \right] du} \nonumber \\
= \frac{1}{1-(1-s)\,t} \, \left[ -s \, \frac{\ln{(1-s)}}{s} +(1-t) \, \frac{\ln{t}}{t-1} \right] \, = \, - \, \frac{\,\ln{[(1-s)\,t]}\,}{\:1-(1-s)\,t\,} \: . \quad
\end{eqnarray}
\begin{flushright} $\Box$ \end{flushright}
\end{prova}

\vspace{0.5cm}

\subsection{Second main result}

\begin{teo}[$\zeta{(3)} \not\in \mathbb{Q}$]
\label{teo:mainz3}
\; ~ The number $\,\zeta{(3)}\,$ is irrational.
\end{teo}

\begin{prova}
\quad  In Lemma~\ref{lem:intPnx}, choose
\begin{equation}
f(x) = - \int_0^1{P_n(y) \: \frac{\,\ln{(x y)}\,}{1 -x y} ~ dy} \, .
\end{equation}
From Lemmas~\ref{lem:linIrrz3} and~\ref{lem:linIrsz3}, for all $\,n \in \mathbb{N}$ we have that
\begin{eqnarray}
J_n = \int_0^1{P_n(x) \, f(x) \: dx} = - \int_0^1\!\!\!\int_0^1{P_n(x) \, P_n(y) \, \frac{\ln{(x y)}}{1-xy} \: dy \, dx} \nonumber \\
\in \: \mathbb{Z} \: \zeta{(3)} + \frac{\mathbb{Z}}{d_n^{\:3}} \, .
\end{eqnarray}
Hence, for some integers $\,a_n$ and $\,b_n$,
\begin{equation}
J_n = \frac{a_n}{d_n^{\:3}} +b_n\,\zeta{(3)} \, .
\label{eq:Inanbnz3}
\end{equation}

On the other hand, from Lemma~\ref{lem:intLn} we have
\begin{eqnarray}
J_n = \int_0^1\!\!\!\int_0^1{P_n(x) \, P_n(y) \: \frac{\left[- \ln{(x y)}\right]}{1-xy} \: dy \, dx} \nonumber \\
= \int_0^1\!\!\!\int_0^1{\,P_n(x) \: P_n(y) \left[ \int_0^1{\frac{1}{1 -(1 -x y)\,z} ~ dz} \right] dy \, dx} \, .
\label{eq:beuk}
\end{eqnarray}
In Ref.~\cite{Beukers79}, Beukers applied repeated integration by parts and a trick change of variables to treat this triple integral. In order to avoid his change of variables, we use the fact that $\,P_n(1-x) = (-1)^n \, P_n(x)\,$ to show that
\begin{equation}
J_n = (-1)^n \, \int_0^1\!\!\!\int_0^1\!\!\!\int_0^1{\,\frac{P_n(x) \, P_n(y)}{1 -[1 -(1-x)\, y]\,z} ~ dx \,dy\,dz} \, .
\end{equation}
From Lemma~\ref{lem:intParcial}, this integral becomes
\begin{eqnarray}
J_n = (-1)^n \, \int_0^1\!\!\!\int_0^1\!\!\!\int_0^1{\,\frac{P_n(x) \, P_n(y)}{[1-(1-z)\,x] \: [1-(1-y)\,z]} ~ dx \,dy\,dz} \nonumber \\
= (-1)^{2n} \, \int_0^1\!\!\!\int_0^1\!\!\!\int_0^1{\,\frac{P_n(x) \, P_n(y)}{[1-(1-z)\,x] \: (1-y z)} ~ dx \,dy\,dz} \, ,
\end{eqnarray}
where we have made use of  $\,P_n(1-y) = (-1)^n \, P_n(y)\,$. From Lemma~\ref{lem:intPnx}, we have
\begin{eqnarray*}
J_n = \int_0^1{\left[\int_0^1{\frac{P_n(x)}{1-(1-z)\,x} ~dx} \: \int_0^1{\frac{P_n(y)}{1-y z} ~ dy} \, \right] dz} \nonumber \\
= \int_0^1 \Bigg[ \frac{(-1)^n}{n!} \, \int_0^1{x^n \, (1-x)^n \frac{\partial^n}{\partial x^n} \left(\frac{1}{1-(1-z)\,x}\right) dx} \; \frac{(-1)^n}{n!} \, \int_0^1{y^n \, (1-y)^n \frac{\partial^n}{\partial y^n} \left(\frac{1}{1-y z}\right) dy} \Bigg] dz .
\end{eqnarray*}
The above partial derivatives are easily determined, yielding
\begin{eqnarray}
J_n = \int_0^1\!\!\!\int_0^1\!\!\!\int_0^1{\,\frac{(x -x^2)^n \, (y -y^2)^n \, (z -z^2)^n}{\left\{[1-(1-z)\,x] \: (1-y z)\right\}^{n+1}} ~ dx \,dy\,dz} \nonumber \\
= \int_0^1\!\!\!\int_0^1\!\!\!\int_0^1{\,\frac{g(x,y,z)^n}{[1-(1-z)\,x] \: (1-y z)} ~ dx \,dy\,dz} ,
\end{eqnarray}
where $\:g(x,y,z) := \dfrac{\:x\,(1 -x) \, y\,(1 -y) \, z\,(1 -z)\:}{[1 -(1 -z)\,x] \; (1 -y z)}$. Clearly, $g(x,y,z)>0\,$ for all inner points of the domain $\,[0,1]^3$, being null only at its boundary.\footnote{Except at the point $(1,1,1)$, where it is an indeterminate form of the kind `$0/0$'. Note, however, that it tends to zero as $(x,y,z) \rightarrow (1^{-},1^{-},1^{-})$.}  Therefore, $|J_n| > 0$, $\forall \, n \in \mathbb{N}$.
On searching for a suitable upper bound, all we need to do is
\begin{eqnarray}
|J_n| = \int_0^1\!\!\!\int_0^1\!\!\!\int_0^1{\,\frac{g(x,y,z)^n}{[1-(1-z)\,x] \: (1-y z)} ~ dx \,dy\,dz} \nonumber \\
\le \int_0^1\!\!\!\int_0^1\!\!\!\int_0^1{\,\max_{[0,1)^3}{\left[ g(x,y,z)^n \right]} \, \frac{1}{[1-(1-z)\,x] \: (1-y z)} ~ dx \,dy\,dz} \nonumber \\
= \left\{ \max_{[0,1)^3}{\left[ g(x,y,z) \right]} \right\}^n \, \int_0^1\!\!\!\int_0^1\!\!\!\int_0^1{\,\frac{1}{[1-(1-z)\,x] \: (1-y z)} ~ dx \,dy\,dz} .
\label{eq:Jn}
\end{eqnarray}
The above maximum can be found analytically by solving the $\,3 \times 3\,$ system
\begin{equation}
\frac{\partial g(x,y,z)}{\partial x} = \frac{\partial g(x,y,z)}{\partial y} = \frac{\partial g(x,y,z)}{\partial z} = 0\, ,
\end{equation}
which reads
\begin{equation}
\left\{
\begin{array}{lllllll}
(1-z)\,x^2 -2 x +1 = 0 \\
z\,y^2 -2\,y +1 = 0 \\
(y-x)\,z^2 -2\,(1-x)\,z +1-x = 0 \, ,
\end{array}
\right.
\end{equation}
respectively. The first two equations are readily solved in terms of $z$, the only solution $\,(x,y) \in (0,1)^2\,$ being
\begin{equation}
\left\{
\begin{array}{lllllll}
x = 1/(1+\sqrt{z}\,) \\
y = \left(1-\sqrt{1-z} \: \right)/z \, ,
\end{array}
\right.
\end{equation}
On substituting these expressions on the third equation, one finds, after some algebra, that
\begin{equation}
\sqrt{z} -z -(\sqrt{z}+z)\,\sqrt{1-z} -z\,\sqrt{z} +1 = 0 \, .
\end{equation}
This simplifies to $\,(\sqrt{z} +z)/\sqrt{1-z} = 1+\sqrt{z}$, which implies that $\,\sqrt{z} =\sqrt{1-z}$, so $\,z = 1/2$. On putting this value of $z$ back in the $(x,y)$ solution, above, one promptly finds $\,x = y = 2 -\sqrt{2}$. So, the point $\,(2 -\sqrt{2},2 -\sqrt{2},1/2)\,$ is the only maximum of $\,g(x,y,z)\,$ in that domain, a result which corrects Eq.~(19) of Ref.~\cite{Miller}, in which the maximum was erroneously found at $\,(2-\sqrt{2}, \sqrt{2}-1, 1/2)$. From Eq.~\eqref{eq:Jn}, one has
\begin{equation}
|J_n| \le \left[g\left(2-\sqrt{2}, 2-\sqrt{2}, \frac12 \right)\right]^n \cdot J_{00} = \left(17 -12\,\sqrt{2}\,\right)^n \times [\,2\;\zeta{(3)}] .
\end{equation}
Since $\,\left(17 -12\,\sqrt{2}\:\right) < 1$, this shows that $\,|J_n| \rightarrow 0\,$ as $n \rightarrow \infty$.  From Eq.~\eqref{eq:Inanbnz3}, one has
\begin{equation}
0 < \left| \frac{a_n}{d_n^{\:3}} + b_n \: \zeta{(3)} \right| \le 2 \left(17 -12\,\sqrt{2}\,\right)^n \zeta{(3)} \, ,
\end{equation}
which means that
\begin{equation}
0 < \left| a_n + b_n \: d_n^{\:3} \: \zeta{(3)} \right| \le \, 2 \: {d_n^{\:3}} \left(17 -12\,\sqrt{2}\,\right)^n \zeta{(3)} \, .
\end{equation}
From Lemma~\ref{lem:dn}, we know that $\,d_n^{\:3} \le (n^{\pi{(n)}})^3 \sim e^{3n} = (e^3)^n$. Since $\,e^3<21$, then, for sufficiently large values of $n$, $d_n^{\:3} < 21^n$, which leads us to
\begin{eqnarray}
0 < \left| a_n + c_n \: \zeta{(3)} \right| < 2 \times 21^n \left(17 -12\,\sqrt{2}\,\right)^n \zeta{(3)}  \nonumber \\
= 2 \left[21 \left(17 -12\,\sqrt{2}\,\right)\right]^n \zeta{(3)} \, ,
\end{eqnarray}
where $\,c_n = b_n \: d_n^{\:3}\,$ is an integer.

Now, assume that $\,\zeta{(3)} \in \mathbb{Q}$. Since $\,\zeta{(3)}>0$, then make $\,\zeta{(3)} = p/q\,$, $\,p$ and $q$ being coprime positive integers.  From the above inequality, one has
\begin{equation}
0 < \left| a_n + c_n \: \frac{p}{q} \right| < 2 \left[21 \left(17 -12\,\sqrt{2}\,\right)\right]^n \, \frac{p}{q} \, .
\end{equation}
On searching for a contradiction, note that $\: 21 \left(17 -12\,\sqrt{2}\,\right) = 0.618\ldots < 2/3\,$, so
\begin{equation}
0 < \left| q \, a_n + p \: c_n \right| < 2 \, p \left( \frac{2}{3} \right)^{\!n} .
\end{equation}
Since $\,| q \, a_n + p \: c_n |\,$ is a positive integer, then $\,\left| q \, a_n + p \: c_n \right| \ge 1\,$ for all $\,n \in \mathbb{N}$. However, $\,2 \, p ~ ( 2/3 )^{n}\,$ is less than $\,1\,$ for sufficiently large values of $\,n$, and we have a contradiction.\footnote{More precisely, $\,\left| q \, a_n + p \: c_n \right| < 1\:$ for all $\:n > \ln{(2\,p)}\,/\ln{(\frac32)} $.}  Therefore, $\zeta{(3)}$ cannot be a positive rational number.
\begin{flushright} $\Box$ \end{flushright}
\end{prova}





\end{document}